\documentclass[10pt]{amsart}

\usepackage{booktabs,setspace}
\usepackage{amsmath,mathtools,amsthm,amsfonts,amssymb,color,verbatim,graphicx}

\usepackage[letterpaper]{geometry}
\usepackage{stmaryrd}

\theoremstyle{definition}

\numberwithin{equation}{section}
\newcommand{\ZZ}{\mathbb{Z}}

\keywords{perfect order subset groups, POS-group, order subsets}

\subjclass[2010]{20D99}

\begin{document}
\setlength{\jot}{0pt} 
\title{On two questions by Finch and Jones about Perfect Order Subset Groups}

\today
\author{Bret Benesh}

\address{
Department of Mathematics,
College of Saint Benedict and Saint John's University,
37 College Avenue South,
Saint Joseph, MN 56374-5011, USA
}
\email{bbenesh@csbsju.edu}

\begin{abstract}
A finite group $G$ is said to be a {\it POS-group} if the number of elements of every order occurring in $G$ divides $|G|$.  We answer two questions by Finch and Jones in~\cite{finch2002curious} by providing an infinite family of nonabelian POS-groups with orders not divisible by $3$.   
\end{abstract}

\maketitle

Let $G$ be a finite group, and define the {\it order subset of an element $x \in G$} to be $\{g \in G \mid o(g)=o(x)\}$, where $o(x)$ denotes the order of $x$.  We say that $G$ has {\it perfect order subsets} if the number of elements in every order subset divides $|G|$; in this case, we say that $G$ is a POS-group.  It is easy to see that $\ZZ_2$, $\ZZ_4$, and the symmetric group $S_3$ are POS-groups, whereas $\ZZ_3$, $\ZZ_5$, and $S_4$ are not.

This definition is due to Finch and Jones, who worked with abelian groups~\cite{finch2002curious} and direct products of abelian groups with $S_3$~\cite{finch2003nonabelian}. They provided the following open questions at the end of~\cite{finch2002curious}.
 
\begin{enumerate}
\item Are there nonabelian POS-groups other than the symmetric group $S_3$? 
\item If the order of a POS-group is not a power of $2$, is the order necessarily divisible by $3$?  This is also Conjecture 5.1 from~\cite{finch2003nonabelian}, also by Finch and Jones.  
\end{enumerate}

The answers are ``yes" and ``no," respectively. The first question was answered by Finch and Jones in ~\cite{finch2003nonabelian}, although all groups were direct products of $S_3$ with an abelian group.   Das~\cite{das2009finite} answered both questions by proving that there exists an action $\theta$ such that the semidirect product $\ZZ_{p^{k}} \rtimes_{\theta} \ZZ_{2^{l}}$ is a POS-group, where $p$ is a Fermat prime, $k \geq 1$, and $2^{l} \geq p-1$.  Feit also answered both questions by indicating that a Frobenius group of order $p(p-1)$ for a prime $p>3$ is a POS-group~\cite{finch2003nonabelian}.  

We now provide an infinite family of groups that simultaneously answers both questions.  Let $n \geq 1$, and consider the group $\ZZ_4 \rtimes \ZZ_{2 \cdot 5^{n}}$ with the inversion action.  The order of this group is $2^3 \cdot 5^{n}$, which is not divisible by $3$.  Consequently, $S_3$ cannot appear as a subgroup or quotient of any of these groups, as $3$ divides the order of $S_3$.

Table~\ref{tab:OrderTable} summarizes the easy calculations required to find the size of each order subset of $\ZZ_4 \rtimes \ZZ_{2 \cdot 5^{n}}$ (one can use geometric sums to verify that all elements of the group are accounted for).  Note that the number of elements of each order divides the order of the group, thereby proving that the groups are POS-groups.    This confirms that $\ZZ_4 \rtimes \ZZ_{2 \cdot 5^{n}}$ is a POS-group. 

\vspace{0.1in}

\begin{table}[h]
\begin{tabular}{lllllll}\toprule
\multicolumn{7}{l}{$\ZZ_{4} \rtimes \ZZ_{2\cdot 5^{n}}$, $|G|=2^3 \cdot 5^{n}$} \\
\midrule
Order: & $1$ & $2$ & $4$ & $5^{m}$ & $2 \cdot 5^{m}$ & $4 \cdot 5^{m}$   \\
Elements: & $1$ & $5$ & $2$    & $4 \cdot 5^{m-1}$ &$4 \cdot 5^{m}$ &$8 \cdot 5^{m-1}$   \\  %
\hline
\end{tabular}
\caption{\label{tab:OrderTable}Here, $n$ is defined so that $n \geq 1$ and $m$ is defined so that $1 \leq m \leq n$.}
\end{table}

\bibliographystyle{amsplain}
\bibliography{MasterBibliography}

\end{document}